\RequirePackage{ifpdf}
\ifpdf % We are running pdfTeX in pdf mode
\documentclass[pdftex]{sigma}
\else
\documentclass{sigma}
\fi

\begin{document}

\renewcommand{\PaperNumber}{037}

\FirstPageHeading

\renewcommand{\thefootnote}{$\star$}

\ShortArticleName{The Research of Thomas P.~Branson}

\ArticleName{The Research of Thomas P.~Branson\footnote{This paper is a
contribution to the Proceedings of the 2007 Midwest
Geometry Conference in honor of Thomas~P.\ Branson. The full collection is available at
\href{http://www.emis.de/journals/SIGMA/MGC2007.html}{http://www.emis.de/journals/SIGMA/MGC2007.html}}}

\Author{Michael G. EASTWOOD~$^\dag$ and A. Rod GOVER~$^\ddag$}

\AuthorNameForHeading{M.G. Eastwood and A.R. Gover}

\Address{$^\dag$~Department of Mathematics, University of Adelaide,
SA 5005, Australia}
\EmailD{\href{mailto:meastwoo@member.ams.org}{meastwoo@member.ams.org}}
\URLaddressD{\url{http://www.maths.adelaide.edu.au/pure/staff/meastwood.html}}

\Address{$^\ddag$~Department of Mathematics,
  The University of Auckland,\\
$\phantom{^\ddag}$~Private Bag 92019,
  Auckland 1,
  New Zealand}
\EmailD{\href{mailto:gover@math.auckland.ac.nz}{gover@math.auckland.ac.nz}}
\URLaddressD{\url{http://www.math.auckland.ac.nz/~gover/}}

\ArticleDates{Received March 27, 2008; Published online April 02, 2008}

\Abstract{The Midwest Geometry Conference 2007 was devoted to the
substantial mathematical legacy of Thomas P. Branson who passed away
unexpectedly the previous year. This contribution to the Proceedings
brief\/ly introduces this legacy. We also take the opportunity of
recording his bibliography. Thomas Branson was on the Editorial Board
of SIGMA and we are pleased that SIGMA is able to publish the
Proceedings.}


\Classification{01A70; 22E46; 22E70;
53A30; 53A55; 53C21;
58J52; 58J60; 58J70}

\bigskip

\rightline{\it In memory of Thomas P. Branson (1953--2006)}

\bigskip

\renewcommand{\thefootnote}{\arabic{footnote}}
\setcounter{footnote}{0}

%\section*{}

Tom's research interests ranged broadly; this is perhaps best
indicated by the diverse group of collaborators listed implicitly in
his bibliography, which follows this section. Tom was also known for
his general knowledge in many areas of mathematics, and many of his
colleagues will remember well his prompt and in-depth replies to email
questions.

In the 1980s Tom's research on conformal invariance was signif\/icantly
ahead of its time. An~enduring theme of his research was the natural
interplay between invariance and the underlying symmetry groups. Tom's
work continues to motivate and inspire a thriving and impressive
international research ef\/fort. Many of the conference speakers have
contributed to the various research trends that Tom Branson started.

It is impossible in a written summary to do justice to Tom's research
career. Here we shall outline just a few directions that we feel to be
especially signif\/icant.

\section[$Q$-curvature and extremal problems]{$\boldsymbol{Q}$-curvature and extremal problems}

Tom Branson is
perhaps most well known for his def\/inition, development, and
application of a~new curvature quantity in Riemannian geometry. For
dimension 4 this f\/irst entered the public arena in his joint work
\cite{22} with {\O}rsted, but it was extended to all even dimensions
and developed signif\/icantly in \cite{28} and \cite{35}; these days it
is usually termed ``Branson's $Q$-curvature''. A key feature is its
conformal transformation law in dimension $n$,
\[e^{-n\omega}\widehat{Q}=Q+P\omega,\]
where $\widehat{g}=e^{2\omega}g$ for a smooth function $\omega$ and
$P$ is the celebrated conformal (i.e.\ conformally covariant) operator
of the form $\Delta^{n/2}+\mbox{\it lower order terms}\,$ due to
Graham--Jenne--Mason--Sparling (and Paneitz in dimension~4). Because $P$
turns out to be a divergence, this generalises the transformation of
Gauss curvature in 2 dimensions and shows that on closed manifolds
$\int Q$ is a~global conformal invariant.

The article with {\O}rsted was motivated by the study of the
functional determinants of integral powers of conformal operators and,
in particular, the issue of extremising such quantities within the
class of conformally related metrics of a f\/ixed volume. This was part
of a theme, extended in his work \cite{26} with Chang and Yang, in
\cite{30} with Gilkey, and in \cite{37}, of developing results for
dimension 4 and higher even dimensions, which paralleled the results
from dimension 2 due to Polyakov and others. An excellent expository
account of these directions, written by Tom Branson himself, is to be
found in \cite{77}, which appears in this volume.

\section[Conformal differential operators]{Conformal dif\/ferential operators}

Tom was also a pioneer
of the systematic construction and study of conformal dif\/ferential
ope\-rators and related issues. The importance and use of ellipticity
was an enduring theme and, for example in~\cite{45}, he identif\/ied the
elliptic operators within the class of second order formally
self-adjoint operators arising as Stein--Weiss gradients. On the other
hand, in \cite{47} he provided a classif\/ication of second order linear
conformal dif\/ferential operators. The work \cite{13}, which constructs
operators between dif\/ferential forms on Minkowski space, was motivated
by their r\^ole as representation intertwinors and links to physics.
This partly followed some earlier work on dif\/ferential forms in
\cite{11}. There he constructed, for example, new conformal
dif\/ferential operators of order 4 and 6 and some applications of these
to variational problems. Notable is that on (unweighted) dif\/ferential
forms the operators he found took the form $L=dMd$, where $d$ is the
exterior derivative, which is itself conformal. Although it was not
highlighted at the time it later became clear that, while conformal
operators on forms with the symbol of $L$ were to be expected,
factorisations along these lines are rare, surprising, and valuable.
This thread was picked up much later in the joint work \cite{65} with
Gover. There, using the Fef\/ferman--Graham ambient metric and its links
to conformal tractor calculus, it was found that on any $k$-forms for
$k\leq n/2-1$ there are conformal dif\/ferential operators generalising
his earlier discoveries, but especially signif\/icant is what is
captured in the details of the factorisation. The operators may be
expressed in the form
\begin{equation}\label{rarefac}
L_k=\underbrace{d*\Bigl\{\, \overbrace{(d*d*)^{n/2-k-1}
+\mbox{\it lower order terms}}^{Q_{k+1}}
\,\Bigr\}}_{G_{k+1}}\,d,
\end{equation}
where $*$ is the Hodge $*$-operator. As an operator on closed forms,
$Q_{k+1}$ generalises the $Q$ curvature and the composition
$d*Q_{k+1}=G_{k+1}$ forms a conformal gauge companion for $d$ or
alternatively for $L_{k+1}$; paired with either of these $G_{k+1}$
yields an injectively elliptic conformal system. Another consequence
of the factorisations (\ref{rarefac}) is that they lead to new
elliptic conformal complexes. A driving motivation here is that such
complexes admit torsions (or determinants), which generalise Cheeger's
de Rham half torsion and, in particular, have Polyakov type conformal
variation formulae. These ideas are sketched in Tom's f\/ine survey
article~\cite{64}.

\section{Spectrum generating functions and intertwining operators}

Another enduring theme of Tom's work was exploiting the implications
of invariance for spectral data. For example, on the round sphere the
conformal covariance relation satisf\/ied by the conformal Laplacian
$D=\Delta+\frac{n(n-2)}{4}$ imposes relations among the eigenvalues of
the Laplacian itself. In fact, in this case there are suf\/f\/icient
relations to deduce all the spectral data, and hence the term
``spontaneous generation of eigenvalues'' which titles the joint work
of Branson with {\O}rsted~\cite{70}. Here, a heavy use is made of
symmetry. The conformal covariance implies a simple formula for the
commutator of $D$ with the Lie derivative along conformal vector
f\/ields. But the sphere has a maximal dimension space of such f\/ields
and, via the commutation relation, one can show that the action of
appropriately chosen conformal vector f\/ields will shift eigenvectors
to linear combinations of eigenvectors with adjacent eigenvalues.

An observation made much earlier by Tom is that, given the explicit
spectral data for a basic Laplacian operator and again suf\/f\/icient
symmetry, one may deduce rather explicit formulae for a series of
intertwinors, which are not generally dif\/ferential. Of course, this
uses that the span of the eigenfunctions for the given Laplacian is
dense on the compact manifold concerned. Once again using the sphere
as an example, one may suppose that some operator $A_{2r}$, acting
between sections of the trivial bundle, is a function of $\Delta$ and
satisf\/ies an obvious generalisation of the conformal covariance
relation enjoyed by $D$. Compressing to eigenspaces of the Laplacian,
this implies characterising relations for the spectrum of the
operator. For this series of operators, Tom was able to obtain the
rather striking formula
\[
A_{2r}=\frac{\Gamma(A_1+\frac{1}{2}+r)}{\Gamma(A_1+\frac{1}{2}-r)},
\qquad A_1=\sqrt{\Delta +\left(\frac{n-1}{2}\right)^2}
\]
for $r\in \mathbb{C}$ and $r\notin \{-n/2,-n/2-1, \dots \}$. These
are intertwinors of the spherical principal series. As pointed out in
\cite{13}, for $r\in \mathbb{Z}_+$ one obtains an explicit formula for
a class of conformal dif\/ferential operators, namely the conformal
Laplacian operators on the sphere:
\[
\prod^r_{p=1} \left\{\Delta +\left(\frac{n}{2}+p-1\right)\left(\frac{n}{2}-p\right) \right\},
\]
(cf.\ Robin Graham's article in this volume\footnote{Graham C.R., Conformal powers of the Laplacian via stereographic projection,
{\it SIGMA} {\bf 3} (2007), 121, 4~pages, \href{http://arxiv.org/abs/0711.4798}{arXiv:0711.4798}.}). In fact, in \cite{13}
Branson used similar ideas to construct the conformal dif\/ferential
operators on forms mentioned earlier. These are dif\/ferential operators
$D_{2l,k}$, on forms of all even orders $2l$, on dif\/ferential forms of
all orders $k$, on the double cover of the $n$-dimensional
compactif\/ied Minkowski space.

Tom's early ideas were put into a representation theoretic framework
and generalised in the signif\/icant work \cite{36}, with \'Olafsson and
{\O}rsted. There they used a spectrum generating operator, constructed
from a combination of quadratic Casimirs and the development of ideas
as above to construct intertwinors for representations induced from a
maximal parabolic subgroup.

\section{Heat kernel asymptotics}

The so-called heat kernel
asymptotics have long been of considerable interest in the study of an
elliptic operator $F$. The coef\/f\/icients of the small time expansion of
such a kernel are spectral invariants; they encode information about
the asymptotic properties of the spectrum of $F$.

Tom made substantial contributions in calculating such invariants.
Perhaps his best known article in this direction is \cite{17} written
joint with Gilkey. They consider a rather general setting, namely a
Riemannian manifold with possible boundary together with a smooth
vector bundle~$V$ with connection~$\nabla$, a given section $E$ of
${\mathrm{End}}(V)$, and an auxiliary smooth function~$f$. The
associated heat kernel is ${\mathrm{Tr}}(fe^{-tP})$, where
$P=\nabla^*\nabla+E$. The inclusion of the function $f$ allows for an
ef\/f\/icient and secure calculation of the coef\/f\/icients $a_n(f,P)$ for
$n=0,\,0.5,\,1,\,1.5,\, 2$ under both Dirichlet and Neumann boundary conditions.
In particular, by setting $f=1$ in $a_2(f,P)$, this f\/ifth coef\/f\/icient
of the heat kernel aysymptotics is calculated for the f\/irst time.
Another appealing feature of incorporating the function $f$ is that
its normal derivatives showing up in their formulae better record the
distribution-like behaviour of the heat kernel near the boundary.

Another pioneering contribution of Tom's work in this area concerned
the heat kernel asymptotics of operators not necessarily of Laplace
type. This joint work with Gilkey and Fulling~\cite{22a}, and later
Avramidi~\cite{55}, assumes the second order operators concerned still
have positive def\/inite leading symbol but does not assume this
coincides with (the inverse of) a background Riemannian metric.
Results, in this general setting, include explicit formulae for the
f\/irst two terms of the asymptotics.

\pdfbookmark[1]{Publications of Thomas P. Branson}{ref}
\renewcommand{\refname}{Publications of Thomas P. Branson}

\LastPageEnding

\end{document}